\documentclass[aos,seceqn,citesort,dvips]{arximspdf}
\usepackage{mathrsfs}

\doi{10.1214/10-AOS811}
\volume{38}
\issue{6}
\pubyear{2010}
\firstpage{3300}
\lastpage{3320}

\makeatletter

\newcommand{\PP}{\mathbb{P}}
\newcommand{\EE}{{\mathbb{E}}}
\newcommand{\eps}{\varepsilon}
\newcommand{\GG}{\mathscr{G}}
\newcommand{\RR}{\mathbb{R}}
\newcommand{\NN}{\mathbb{N}}
\newcommand{\sub}{\subseteq}
\newcommand{\XX}{\mathscr{X}}
\newcommand{\PPP}{\mathscr{P}}
\newcommand{\HHH}{\mathbb{H}}
\newcommand{\BBB}{\mathbb{B}}
\newcommand{\CCC}{\mathbb{C}}
\renewcommand{\phi}{\varphi}

\newcommand{\given}{\mid}

\newproclaim{df}{Definition}[section]
\newtheorem{theorem}[df]{Theorem}
\newtheorem{cor}[df]{Corollary}
\newtheorem{lem}[df]{Lemma}

\makeatother

\begin{document}
\begin{frontmatter}

\title{Adaptive nonparametric Bayesian inference using location-scale mixture priors}
\runtitle{Adaptive Bayesian inference using mixtures}

\thankstext{t1}{Supported by the Netherlands Organization for Scientific Research NWO.}

\begin{aug}
\author[A]{\fnms{R.} \snm{de Jonge}\ead[label=e1]{r.d.jonge@tue.nl}} and
\author[A]{\fnms{J. H.} \snm{van Zanten}\corref{}\ead[label=e2]{j.h.v.zanten@tue.nl}\thanksref{t1}}
\runauthor{R. de Jonge and J. H. van Zanten}
\affiliation{Eindhoven University of Technology}
\address[A]{Department of Mathematics\\
P.O. Box 513\\
5600 MB Eindhoven\\
The Netherlands\\
\printead{e1} \\
\phantom{E-mail: }\printead*{e2}} 
\end{aug}

\received{\smonth{7} \syear{2009}}
\revised{\smonth{12} \syear{2009}}

%
\begin{abstract}
We study location-scale mixture priors for nonparametric statistical
problems, including multivariate regression, density estimation and
classification. We show that a rate-adaptive procedure can be obtained
if the prior is properly constructed. In particular, we show that
adaptation is achieved if a kernel mixture prior on a regression
function is constructed using a Gaussian kernel, an inverse gamma
bandwidth, and Gaussian mixing weights.
\end{abstract}

%

\end{frontmatter}

\section{Introduction}

In Bayesian nonparametrics, the use of location-scale mixtures of
kernels for the construction of priors on probability densities
is well esthablished. The methodology is used in a variety of
practical settings, and in recent years there has been substantial
progress on the the mathematical, asymptotic theory for kernel
mixture priors as well; cf.
\cite{GGR99,Tok06,GvdVMixtures,MR2336864,WuGhosal,Willem}.
At the present time, we have a well-developed understanding of important
aspects including consistency, convergence rates, rate-optimality
and adaptation properties. A similar, parallel development has
taken place in the area of beta mixture priors; cf.
\mbox{\cite{PetroneWasserman,GhosalBernstein,KruijerVaart,Rou08}}.

A discrete location-scale mixture of a fixed probability density $p$ on
$\RR^d$ can
be expressed as
%
%
\begin{equation}\label{eq: mix1}
x \mapsto\sum_{j=1}^m w_j \frac1{\sigma^d} p \biggl(\frac
{x-x_j}{\sigma} \biggr),
\end{equation}
where $m \in\NN$, $x_1, \ldots, x_m \in\RR^d$, $w_1, \ldots, w_m
\ge0$ and $\sum w_j = 1$,
and $\sigma> 0$.
A prior on densities is obtained by
putting prior distributions on $m$, the locations $x_j$, the scale
$\sigma$ and the
weights $w_j$.
When $p$ satisfies some regularity conditions, a wide class of probability
densities can be well approximated by mixtures of the form (\ref{eq: mix1}).
This indicates that if the priors on the coefficients are
suitably chosen, the resulting prior and posterior on probability
densities can be expected
to have good asymptotic properties. The cited papers
give precise conditions under which this is indeed the case.

Obviously, a much wider class of functions is well approximated by mixtures
of the form (\ref{eq: mix1})
if we lift the restriction that the weights $w_j$ should be nonnegative
and sum up to $1$.
This suggests that location-scale mixtures might be attractive priors
not just in the
setting of density estimation, but for instance also in nonparametric
regression.
Although this idea has been proposed in the applied literature; cf.,
for example,~\cite{Short,Higdon}, it does not seem to
have attracted a great deal of attention.
The few examples do show however that the approach can yield quite
satisfactory results.

In the paper~\cite{Short}, location-scale mixture priors
are used in an astrophysical setting for the analysis of data from
galatic radio sources.
The statistical problem essentially boils down to a bivariate,
nonparametric, fixed design regression problem.
The use of a mixture prior is natural in that particular application because
it reflects the idea that the function of interest, which describes the
strength of the
magnetic field caused by our planet and its ``neighborhood'' in space,
is in fact an aggregate of contributions from a large number of
locations, with different weights, which
can be positive
or negative.

Another reason for using a location-scale mixture prior in multivariate
regression, instead of
for instance the popular Gaussian squared exponential or Mat\'ern
priors, are computational advantages.
Conditional on the gridsize $m$ the prior only involves finitely many
terms, so no artificial truncation or approximation is necessary for
computation.
As argued also in
\cite{Short}, the mixture prior allows to avoid the inversion or
decomposition of nontrivial and often ill-behaved
$n \times n$ matrices (with $n$ the sample size),
which can become cumbersome already for moderate sample sizes (cf.
also the discussion in~\cite{Banerjee08}).
In the astrophysical application of~\cite{Short},
the sample size is of the order $1500$ and it is shown that samples of
this order can be dealt
with effectively using kernel mixture priors.

On the theoretical side, little or nothing seems to be known for kernel
mixture priors
in a regression setting.
In the present paper, we therefore take up the study of
asymptotic properties, in order to assess the fundamental
potential of the methodology and to provide a theoretical underpinning
of its use in practice.
We will show that if the kernel and the
priors on locations and scales are appropriately chosen,
kernel mixture priors yield posteriors with very good asymptotic properties.
It is well known that for the estimation of an $\alpha$-regular
function of $d$ variables,
the best possible rate of convergence is of the order $n^{-\alpha
/(d+2\alpha)}$, where
$n$ is the number of observations available.
We will prove that up
to a logarithmic factor, this optimal rate can be attained with
location-scale mixture priors.
More importantly, the near optimal rate can be achieved by
a prior that does \textit{not} depend on the unknown smoothness level
$\alpha$ of the regression function.
In other words, we can obtain a fully adaptive procedure.

The bounds for the convergence rates that we will obtain depend
crucially on the smoothness
of the kernel $p$ that is used. For kernels with only a finite degree
of regularity,
we get suboptimal rates. We only obtain the optimal minimax rate (up to
a logarithmic factor)
for kernels that are infinitely smooth, in the sense that they
admit an analytic extension to a strip in complex space. The standard
normal kernel is an example
of an optimal choice in this respect. We also have to put (mild)
conditions on the priors on the
grid size $m$ and the scale $\sigma$. In particular, the popular
inverse gamma choice for the
scale is included in our setup.

Perhaps surprising is the fact that although we use a probability
density $p$
to construct the mixtures, we can still achieve adaptation to all
smoothness levels.
Intuition from kernel estimation might suggest that when $p$ is a
centered probability density, we have good approximation behavior for
regression functions
with regularity at most $2$, and that for more regular functions we
should use
higher order kernels. This turns out not to be the case however. To
prove this fact,
we adapt an observation of Rousseau, who uses a similar idea to prove
that for densities on the unit interval, using appropriate mixtures of
beta densities
yields adaptation to all smoothness levels; see~\cite{Rou08}.
The recent preprint~\cite{Willem}, which was written at the same time
and independently
of the present work, employs the same idea
to prove adaptation for kernel mixture priors for density estimation.
In the present paper, we extend the technique to a multivariate setting
(see Lemma~\ref{lem: judith} ahead).

The literature on Bayesian
adaptation is still relatively young.
Earlier papers include~\cite{BelitserGhosal,Huang,GLvdV08,GLvdV,LembervdV,Rou08}
and~\cite{VaartZanten09}.
Priors that yield adaptation across a continuum of regularities
in nonparametric regression have been exhibited in~\cite{Huang},
where priors based on spline expansions are considered, and
\cite{VaartZanten09},
which uses randomly rescaled Gaussian processes as priors.

The location-scale priors we consider in this paper
are conditionally Gaussian, since we will put Gaussian priors on the
mixing weights.
This allows us to use the machinery for Gaussian process priors
developed in
\cite{vdVvZ} and~\cite{vdVvZRKHS} in our proofs. Other technical
ingredients include metric entropy results for
spaces of analytic functions, as can be found, for instance, in
\cite{Kolmogorov},
and the connection between metric entropy and small
deviations results for Gaussian process (cf.
\cite{KuelbsLi,LiLinde}).
We will obtain a general result for a conditionally Gaussian kernel
mixture process, which can in fact be used in a variety
of statistical settings.
To illustrate this, we present rate of contraction results not just for
nonparametric regression,
which is our main motivation, but also for density estimation and
classification settings.

In the next section, we present the main results of the paper. In
Section~\ref{sec: general},
we state a general result for a conditionally Gaussian location-scale
mixture process whose
law will be used to define the kernel mixture prior in the various
statistical settings. Rate of contraction results for
nonparametric regression, density estimation and classification are
given in Section~\ref{sec: statistics}.
The proof of the general theorem can be found in Section~\ref{sec: proof}.

\subsection{Notation}

\begin{itemize}
\item
$\Im z$, $\Re z$: imaginary and real part of a complex number $z$.
\item
$\NN_0 = \NN\cup\{0\}$.
\item
For $k \in\NN_0^d$: $k. = k_1+\cdots+ k_d$, $k! = k_1! \cdots k_d!$.
\item
$f*g$: convolution of $f$ and $g$.
\item
$a \vee b = \max\{a,b\}$, $a\wedge b = \min\{a,b\}$, $a_+ = a\vee0$.
\item
$C(X)$: continuous functions on $X$.
\item
$C^\alpha(X)$ for $\alpha> 0$ and $X \sub\RR^d$: functions on $X$
with bounded partial derivatives
up to the order $\beta$, which is the largest integer strictly smaller
than $\alpha$, and such
that the partial derivatives of order $\beta$ are H\"older continuous
of order $\alpha-\beta$.
For $f \in C^\alpha(X)$ we denote by $\|f\|_\alpha$ the associated H\"
older norm of $f$;
cf.~\cite{vdVW}, Section 2.7.1. The H\"older ball of radius $R > 0$
is defined as $C^\alpha_R(X) =
\{f \in C^\alpha(X)\dvtx\|f\|_\alpha\le R\}$.
\end{itemize}

\section{Main results}

\subsection{General result for Gaussian location-scale mixtures}
\label{sec: general}

On a common probability space, let $M$ be an $\NN$-valued random variable,
$\Sigma$ a $(0,\infty)$-valued random variable and
$(Z_k\dvtx k \in\NN^d)$ standard Gaussian random variables, all independent.
The stochastic process $W$ indexed by $[0,1]^d$ is defined by
%
%
\begin{equation}\label{eq: rf}
W(x) = \sum_{k \in\{1, \ldots, M\}^d} Z_k\frac{1}{M^{d/2}}\frac
1{\Sigma^d}p \biggl(\frac{x-k/M}{\Sigma} \biggr)
\end{equation}
for $x \in[0,1]^d$, where $p\dvtx\RR^d \to\RR$ is a function that
belongs to the class $\PPP_\gamma$ of
$\gamma$-regular kernels defined as follows.
\begin{df}\label{df: p}
For $\gamma\in(d/2, \infty]$, an integrable function $p$ on $\RR^d$
belongs to $\PPP_\gamma$ if
$\int_{\RR^d} p(x) \,dx = 1$, it is uniformly Lipschitz on $\RR^d$,
it has finite moments of every order,
and it satisfies one of the following conditions, depending on whether
$\gamma< \infty$ or $\gamma= \infty$:
\begin{itemize}
\item
For $\gamma< \infty$: $p$ belongs to $C^\gamma(\RR^d)$.
\item
For $\gamma= \infty$: $p$ is the restriction to $\RR^d$ of a
function that is defined on the
set $S = \{(z_1, \ldots, z_d) \in\CCC^d\dvtx|\Im z_{j}|\le1$ for $j =
1, \ldots, d\}$,
and that is bounded and analytic on $S$.
\end{itemize}
\end{df}

Examples of kernels belonging to $\PPP_\gamma$ for $\gamma< \infty$
are abundant.
Using Fourier inversion, it is not difficult to see that an integrable
function $p$ belongs to $\PPP_\infty$ if it has a
characteristic function
\[
\psi(\lambda) = \int_{\RR^d} e^{i(\lambda, x)}p(x) \,dx,
\]
which is infinitely often differentiable at $0$, which satisfies $\psi
(0) = 1$, and which satisfies
the exponential moment condition
\[
\int_{\RR^d} e^{\|\lambda\|} |\psi(\lambda)| \,d\lambda< \infty.
\]
The prime example is the standard normal density on $\RR^d$, which is
easily seen to belong to $\PPP_\infty$.
Note that we do not require that
$p \ge0$ in Definition~\ref{df: p}. So, in fact, higher order kernels
are allowed as well.

The index $\gamma$ of the class of kernels quantifies the regularity
of the kernel that is employed.
We will see that this regularity influences the rate of convergence
that we can obtain for the corresponding
location-scale mixture prior.
The restriction $\gamma> d/2$ is connected to the fact that in order
to obtain bounds for the process
$W$ independent of $M$, we want the process in (\ref{eq: rf}) to be
well defined if the sum is
taken over all $k$ in $\NN^d$.

For $\eps>0$, the metric entropy of a set $B$ in a metric space with
metric $d$ is defined as $\log N(\eps,B,d)$, where $N(\eps,B,d)$ is
the minimum number of balls of radius $\eps$ needed to cover $B$.
Fix $0 < a < b < 1$ and define $\XX= [a,b]^d$. Let $d_\gamma=
2d(d+\gamma)/(2\gamma- d)$ and
$\delta_\gamma= d/(2\gamma- d)$.
\begin{theorem}\label{thm: general}
Suppose that $p \in\PPP_\gamma$ for $\gamma\in(d/2,\infty]$, that
$\PP(M = m) \ge Cm^{-s}$ for some $C> 0$, $s > 1$,
and that $\Sigma$ has a Lebesgue density $g$ that, for some
$D_{1},D_{2},D_{3},D_{4}>0$ and
$q,r \ge0$, satisfies
%
%
\begin{equation}\label{eq: g}
D_{1}\sigma^{-q}e^{-D_{2}({1/\sigma})^{d_\gamma}(\log
{1/\sigma})^{r}} \le g(\sigma)
\le D_{3}\sigma^{-q}e^{-D_{4}({1/\sigma})^{d_\gamma}(\log
{1/\sigma})^{r}}
\end{equation}
for all $\sigma$ in a neighborhood of $0$.

Then if $w_0 \in C^\alpha(\XX)$ for $\alpha> 0$, there exist
for every constant $C > 1 $ measurable subsets $B_n$ of $C([0,1]^d)$
and a constant $D > 0$ such that, for $n$ large enough,
%
%
\begin{eqnarray}
\label{eq: toshow1}\log N(\overline\eps_n, B_n, \|\cdot\|_\infty)
& \le & D n\overline\eps^2_n,\\
\label{eq: toshow2}\PP(W \notin B_n) & \le & e^{-Cn\eps^2_n},\\
\label{eq: toshow3}\PP\Bigl( {\sup_{x \in\XX}}|W(x)-w_0(x)| \le\eps
_n \Bigr) & \ge & e^{-n\eps^2_n}.
\end{eqnarray}
Here if $\gamma< \infty$,
\[
\eps_n = {n}^{-{\alpha}/({d_\gamma+2\alpha(1+\delta_\gamma)})},\qquad
\overline\eps_n = n^{-({\alpha(1-(d\delta_\gamma)/(2\gamma))})/
({(d_\gamma+ 2\alpha(1+\delta_\gamma))(1+d/(2\gamma))})},
\]
and if $\gamma= \infty$,
\begin{eqnarray*}
\eps_n &=& {n}^{-{\alpha}/({d+2\alpha})} \log^{({r\vee
(1+d)})/({2+d/\alpha})} n,\\
\overline\eps_n &=& n^{-{\alpha}/({d+2\alpha})}\log^{
({r\vee(1+d)})/({2+d/\alpha})
+({1+d-r})/{2}_+} n.
\end{eqnarray*}
\end{theorem}

A few remarks about
the result are in order.
First of all, the process $W$ is indexed by the unit cube, but the
supremum in
(\ref{eq: toshow3}) is over the strictly smaller set $\XX$. This is
due to the fact that to obtain
good enough approximations of the given function $w_0$ defined on $\XX$
by location-scale mixtures of the kernel $p$, we also need kernels
centered at points
just outside the set $\XX$. A result like (\ref{eq: toshow3}) with
the supremum over the entire
unit cube is only possible under additional assumptions on the boundary behavior
of the function~$w_0$.

Theorem~\ref{thm: general} connects to existing results
for nonparametric Bayes procedures, which give
sufficient conditions of the form (\ref{eq: toshow1})--(\ref{eq:
toshow3}) for
having a certain rate of posterior contraction; cf., for example,
\cite{GGvdV,GvdVnoniid,Meulen06}. In the next subsection, we will
single out the most important
particular cases. In all cases, the statistical results will state that
the posterior will
asymptotically concentrate on balls of radius of the order
$\overline\eps_n$ around the true parameter (relative to a natural
statistical metric
depending on the specific setting). Note that in the case $\gamma<
\infty$, this means we
only obtain a rate if $(d\delta_\gamma)/(2\gamma) < 1$, which is
true if
and only if $\gamma> (1/4)(1+\sqrt{5})d \approx(0.81) d$.
In particular, the choice $\gamma\ge d$ suffices to have consistency.
As the smoothness $\gamma$ of the kernel $p$ that is employed is increased,
the rate of contraction improves. Since $d_\gamma\to d$ and $\delta
_\gamma\to0$
as $\gamma\to\infty$, the power of $n^{-1}$ in the expression for
the rate $\overline\eps_n$
tends to $\alpha/(d+2\alpha)$ as $\gamma\to\infty$, which
corresponds to
the optimal minimax rate of convergence for estimating an $\alpha
$-regular function of $d$ variables.
If an analytic kernel $p \in\PPP_\infty$ is used the minimax rate
$n^{-\alpha/(d+2\alpha)}$ itself is attained, up
to a logarithmic factor.

The proof of the theorem is deferred to Section~\ref{sec: proof}. In
the next subsection, we give the precise rate of contraction result for
nonparametric regression, density estimation and classification settings.
The first case, which was the original motivation for this study, is
worked out in some detail.
The analogous results for the second and third settings are presented
more briefly, to avoid unnecessary duplications.

\subsection{Rate of contraction results for specific statistical settings}
\label{sec: statistics}

\subsubsection{Regression with Gaussian errors}

Consider a multivariate regression problem where we have
known design points $x_1, x_2, \ldots\in\XX=[a,b]^d$ for some $a<
b$ and $d \in\NN$, and we observe
real-valued variables $Y_1, \ldots, Y_n$ satisfying the regression relation
\[
Y_i = \theta(x_i) + \eps_i
\]
for $\theta\dvtx\XX\to\RR$ an unknown regression function and error variables
$\eps_i$ that are independent and Gaussian, with mean $0$ and variance
$\tau^2$.
We assume that $0 < a< b< 1$, so that the design space $\XX$ is
strictly contained in the
interior of the unit cube in $\RR^d$.

As prior on the regression function, we employ the law $\Pi_\Theta$
that the stochastic process
$W$ defined by (\ref{eq: rf}) generates on the space $C(\XX)$ of
continuous functions on $\XX$.
The total prior $\Pi$ on the pair
$(\theta, \tau)$ is then defined by $\Pi(d\theta, d\tau) = \Pi
_\Theta(d\theta) \times\Pi_T(d\tau)$,
for $\Pi_T$ a prior on a compact interval that is assumed to contain
the true value $\tau_0$,
with a Lebesgue density
that is bounded away from $0$.

The posterior distribution
for $(\theta, \tau)$ given the data $Y_1, \ldots, Y_n$ is denoted by
$\Pi(\cdot\given Y_1, \ldots, Y_n)$. By Bayes formula, it is given
by the expression
\[
\Pi(B \given Y_1, \ldots, Y_n) = \frac{\int_B L(\theta, \tau;
Y_1, \ldots, Y_n) \Pi(d\theta, d\tau)}
{\int L(\theta, \tau; Y_1, \ldots, Y_n) \Pi(d\theta, d\tau)},
\]
where
\[
L(\theta, \tau; Y_1, \ldots, Y_n) = \frac{1}{(2\pi\tau^2)^{n/2}}
\exp\Biggl(-\frac1{2\tau^2}\sum_{i=1}^n\bigl(Y_i-\theta(x_i)\bigr)^2 \Biggr)
\]
is the likelihood.
For a given sequence of positive numbers $\eps_n \downarrow0$,
the posterior is said to contract around the true parameter $(\theta
_0, \tau_0)$
at the rate $\eps_n$ if for $L > 0$ sufficiently large,
\[
\Pi\Biggl((\theta, \tau)\dvtx\frac1n\sum_{j=1}^n\bigl(\theta(x_j)
-\theta_0(x_j)\bigr)^2 + |\tau- \tau_0|^2 > L^2\eps^2_n \given Y_1,
\ldots,
Y_n \Biggr)\stackrel{P_{(\theta_0, \tau_0)}}{\longrightarrow} 0
\]
as $n \to\infty$,
where the convergence is in probability under the true distribution
governed by $(\theta_0, \tau_0)$.
This means in particular that asymptotically, the marginal posterior
for $\theta$ is concentrated on balls with radius of the order $\eps
_n$ around
the true regression function $\theta_0$, where we use the natural $L^2$-norm
associated to the empirical measure of the design points to measure distance.

The next theorem follows from Theorem~\ref{thm: general}, in
combination with
the results in~\cite{GvdVnoniid} (slightly adapted like Theorem 2.1 of
\cite{GvdVMixtures}
in the density estimation
case; cf. also the discussion following Theorem 3.1 of~\cite{VaartZanten09}).
\begin{theorem}\label{thm: reg}
Suppose that the conditions of Theorem~\ref{thm: general} are fulfilled.
Then if $\theta_0 \in C^\alpha(\XX)$ for $\alpha> 0$, the posterior
contracts
at the rate
\[
n^{-{\alpha(1-(d\delta_\gamma)/(2\gamma))}/
({(d_\gamma+ 2\alpha(1+\delta_\gamma))(1+d/(2\gamma))})},
\]
if $\gamma< \infty$, or at the rate
\[
n^{-{\alpha}/({d+2\alpha})}\log^{({r\vee(1+d)})/({2+d/\alpha})
+({1+d-r})/{2}_+} n,
\]
if $\gamma= \infty$.
\end{theorem}

As discussed above already the choice $p \in\PPP_\infty$ yields the
best rate of contraction,
namely the optimal minimax rate, up to a logarithmic factor.
Also note that the prior does not depend on the unknown regularity
$\alpha$ of the true regression function,
so the procedure is rate-adaptive.
Observe that for $p \in\PPP_\infty$ and $r = 1+d$ we obtain the rate
$(n/\log^{1+d}n)^{-\alpha/(d+2\alpha)}$.
If $r$ is strictly larger or smaller than $1+d$, we get a slightly
worse rate, in the
sense that the power of the logarithm in our upper bound for the rate increases.

In the following corollary, we single out the important special case of
a standard Gaussian kernel and
an inverse gamma prior (or a power of it in the multivariate case) on
the scale.
\begin{cor}\label{cor: reg}
Suppose that $p$ is the standard Gaussian
density on $\RR^d$, $\Sigma^d$ is inverse gamma, and $M$ is such that
$\PP(M = m) \ge C m^{-s}$ for some $C > 0$ and $s > 1$.
Then if $\theta_0 \in C^\alpha(\XX)$ for $\alpha> 0$, the posterior
contracts
at the rate
\[
{n}^{-{\alpha}/({d+2\alpha})} \log^{({4\alpha+ 4\alpha d +
d + d^2})/({4\alpha+ 2d})} n.
\]
\end{cor}
\begin{pf}
Simply note that the standard normal kernel belongs to $\PPP_\infty$
and that
if $\Sigma^d$ has an inverse gamma law, then (\ref{eq: g}) is
satisfied with
$r = 0$.
\end{pf}

\subsubsection{Density estimation}

Let $X_1, \ldots, X_n$ be a sample from a positive density $f_0$ on
the set $\XX= [a,b]^d$, for
$0< a< b < 1$. The aim is to estimate the unknown density.

We consider the prior $\Pi$ on densities defined as the law that is
generated on the function space $C(\XX)$
by the random function
%
%
\begin{equation}\label{eq: dens}
x\mapsto\frac{e^{W(x)}}{\int_\XX e^{W(y)}\, dy}
\end{equation}
for $W$ the process defined by (\ref{eq: rf}). In this case, we say
that the posterior
$\Pi(\cdot\given X_1, \ldots, X_n)$ contracts around the true
density $f_0$
at the rate $\eps_n$ if for all $L > 0$ large enough,
\[
\Pi\bigl(f\dvtx h(f,f_0) > L\eps_n \given X_1, \ldots, X_n\bigr) \stackrel
{P_{f_0}}{\to} 0
\]
as $n \to\infty$, where $h$ is the Hellinger distance.

Theorem~\ref{thm: general}, the general rate of contraction results
for Bayesian density estimation (cf.~\cite{GGvdV,GvdVMixtures})
and the relations between the uniform norm on the paths of $W$ and the
relevant statistical metrics
on the densities (\ref{eq: dens}) (cf.~\cite{vdVvZ}) yield the
following result.
\begin{theorem}
In this setting, the assertions of Theorem~\ref{thm: reg} and
Corollary~\ref{cor: reg}
are true for $\theta_0 = \log f_0$.
\end{theorem}

\subsubsection{Classification}

Consider i.i.d. observations $(X_1, Y_1), \ldots, (X_n, Y_n)$, where the
$X_i$ take values in the set $\XX= [a,b]^d$, $0 < a< b< 1$, and
the $Y_i$ take values in $\{0,1\}$. The aim is to estimate the regression
function $r_0(x) = \PP(Y_1 = 1\given X_1 = x)$.

As prior on $r_0$, we use the law $\Pi$ of the process $\Psi(W)$,
where $W$ is
as in~(\ref{eq: rf}) and the link function $\Psi\dvtx\RR\to(0,1)$ is
the logistic or normal distribution
function. Let $\Pi(\cdot\given(X_1, Y_1), \ldots, (X_n, Y_n))$
denote the corresponding posterior
and let $G$ be the distribution of the covariate $X_1$. With $\|\cdot\|
_{2,G}$ the associated $L^2$-norm,
we say that the posterior contracts around the truth $r_0$ at the rate
$\eps_n$ if for all large enough $L> 0$,
\[
\Pi\bigl(r\dvtx\|r-r_0\|_{2,G} > L\eps_n \given(X_1, Y_2), \ldots, (X_n,
Y_n)\bigr) \stackrel{P_{r_0}}{\to} 0
\]
as $n \to\infty$.

Theorem~\ref{thm: general}, the general rate of contraction results
(cf.~\cite{GGvdV})
and the relations between the relevant norms (cf.~\cite{vdVvZ})
yield the following result.
\begin{theorem}
In this setting, the assertions of Theorem~\ref{thm: reg} and
Corollary~\ref{cor: reg}
are true for $\theta_0 = \Psi^{-1}(r_0)$.
\end{theorem}

\section{\texorpdfstring{Proof of Theorem \protect\ref{thm: general}}{Proof of Theorem 2.2}}
\label{sec: proof}

We will find the appropriate sieves $B_n$ and derive the inequalities
(\ref{eq: toshow1})--(\ref{eq: toshow3})
by using the fact that conditionally on the grid size $M$ and the scale
$\Sigma$,
the process $W$ is Gaussian. For fixed $m \in\NN$ and $\sigma> 0$,
we define the stochastic process $(W^{m,\sigma}(x)\dvtx x \in[0,1]^d)$ by
setting
\[
W^{m,\sigma}(x) = \sum_{k \in\{1, \ldots, m\}^d} Z_k\frac
{1}{m^{d/2}}\frac1{\sigma^d}p \biggl(\frac{x-k/m}{\sigma} \biggr).
\]
In the following subsection, we first study some properties of the
Gaussian process $W^{m,\sigma}$
that we will need to establish (\ref{eq: toshow1})--(\ref{eq: toshow3}).

\subsection{\texorpdfstring{Properties of $W^{m,\sigma}$}{Properties of W m, sigma}}

Recall that in general, the reproducing kernel Hilbert space (RKHS)
$\HHH$ attached to a zero-mean Gaussian process $X$
is defined as the completion of the linear space of functions $t\mapsto
\EE X({t})H$ relative to the inner product
\[
\langle\EE X(\cdot)H_{1},\EE X(\cdot)H_{2}\rangle_{\HHH}=\EE H_{1}H_{2},
\]
where $H$, $H_{1}$ and $H_{2}$ are finite linear combinations of the
form $\sum_{i}a_{i}X({s_{i}})$ with $a_{i}\in\RR$ and $s_{i}$ in the
index set of $X$. The following lemma describes the RKHS
of the process $W^{m,\sigma}$. It is a direct consequence of a general
result describing the
RKHS of a Gaussian process admitting a series expansion; cf. Theorem
4.2 of~\cite{vdVvZRKHS} and the
discussion following it.
\begin{lem}
\label{lem: RKHS}
The reproducing kernel Hilbert space $\HHH^{m,\sigma}$ of $
W^{m,\sigma}$ consists of all functions
of the form
%
%
\begin{equation}
\label{rkhs}
h(x)=\sum_{k \in\{1,\ldots,m\}^{d}} w_{k}\frac{1}{\sigma^{d}}
p \biggl(\frac{x-k/m}{\sigma} \biggr),\qquad x\in[0,1]^{d},
\end{equation}
where the weights $w_{k}$ range over the entire set of real numbers.
The RKHS-norm is given by
%
%
\begin{equation}
\label{def rkhs-norm}
\|h\|_{\HHH^{m,\sigma}}^{2}= m^{d} \min_{w} \sum_{k \in\{
1,\ldots,m\}^{d}} w_{k}^{2},
\end{equation}
where the minimum is over all weights $w_k$ for which the
representation (\ref{rkhs}) holds true.
\end{lem}

We remark that if the functions $x \mapsto p((x-k/m)/\sigma)$ on
$[0,1]^d$ are
linearly independent, then the representation (\ref{rkhs}) of an
element of the RKHS
is necessarily unique and hence the minimum in (\ref{def rkhs-norm})
can be removed.
For our purpose, it is, however, not important that these functions are
independent for
every fixed $\sigma$ and~$m$.

Next, we consider the so-called centered small ball probabilities of
the process $W^{m,\sigma}$, which are determined by
its reproducing kernel Hilbert space. We use well-known results by
Kuelbs and Li~\cite{KuelbsLi} and Li and Linde~\cite{LiLinde} that
relate the metric entropy of the unit ball in the RKHS to the centered
small ball probabilities of the process.
The unit ball $\HHH^{m,\sigma}_1$ in the reproducing kernel Hilbert
space $\HHH^{m,\sigma}$ is the set of all elements
$h \in\HHH^{m,\sigma}$ such that $\|h\|_{\HHH^{m,\sigma}}\le1$.

To find an upper bound for the metric entropy of the unit ball, we
embed it in appropriate space of functions for which
an upper bound for the entropy is known, depending on the value of
$\gamma$.
First, we consider the case $\gamma< \infty$.
Let $h$ be an element of $\HHH^{m,\sigma}$. By Lemma~\ref{lem:
RKHS}, it admits a representation~(\ref{rkhs}), with the weights $w_k$ such that $\|h\|^2_{\HHH
^{m,\sigma}} = m^d\sum w^2_k$.
If $p \in\PPP_\gamma$ with $\gamma< \infty$, we get that $h \in
C^\gamma([0,1]^d)$ and
$\|h\|_\gamma\le{\sigma^{-(d+\gamma)}}{\|p\|_\gamma}\|h\|_{\HHH
^{m,\sigma}}$.
Hence, we have $\HHH^{m,\sigma}_1 \subset C^\gamma_{R}([0,1]^d)$ in
this case, where
$R = {\sigma^{-(d+\gamma)}}{\|p\|_\gamma}$.
For $\gamma= \infty$ and $h$ as before, it follows from the
assumptions on $p$ that the function $h$ is in fact well defined on
$S_{\sigma}=\{z\in\CCC^{d}\dvtx\forall j$ $|\Im z_{j}|\le\sigma\}$,
is analytic on this set and takes real values on $\RR^{d}$.
By the Cauchy--Schwarz inequality, it follows that
\[
|h(z)|^{2}\le\frac{1}{\sigma^{2d}} \biggl(\sum_{{k \in\{1, \ldots,
m\}^d}}w_{k}^{2} \biggr)
\biggl(\sum_{k \in\{1, \ldots, m\}^d} \biggl|p \biggl(\frac
{z-k/m}{\sigma} \biggr) \biggr|^{2} \biggr).
\]
The last factor on the right-hand side is bounded from above by a
multiple of $m^{d}$ on the set $S_{\sigma}$. Hence, we obtain
%
%
\begin{equation}
\label{grens voor h}
|h(z)|\le K\sigma^{-d}\|h\|_{\HHH^{m,\sigma}}
\end{equation}
for every $z\in S_{\sigma}$, where the constant $K$ only depends on
the density $p$.
Let $\GG_{\sigma}$ the set of all analytic functions on\vadjust{\goodbreak} $S_{\sigma
}$, uniformly bounded by $K\sigma^{-d}$ on that set,
with $K$ the same constant as in (\ref{grens voor h}).
The preceding shows that for the RKHS unit ball we have $\HHH
^{m,\sigma}_{1} \subset\GG_{\sigma}$ if $\gamma= \infty$.

We see that in all cases we can embed the RKHS unit ball $\HHH
^{m,\sigma}_1$ in a function space
independent of $m$, for which the metric entropy relative to the
supremum norm on $[0,1]^d$ is essentially known. We have the following
result.
\begin{lem}\label{lem: metric entropy}
If $\gamma< \infty$, then
\[
\log N\bigl(\eps,C^\gamma_{{\sigma^{-(d+\gamma)}}{\|p\|_\gamma
}}([0,1]^d),\|\cdot\|_{\infty}\bigr) \le K_0 \biggl(\frac{1}{\eps\sigma
^{d+\gamma}} \biggr)^{{d}/{\gamma}}
\]
for all $\sigma,\eps> 0$, with $K_0$ a constant independent of $\eps
, m$ and $\sigma$.

There exist $\eps_0, \sigma_0 > 0$ such that
\[
\log N(\eps,\GG_{\sigma},\|\cdot\|_{\infty}) \le K_1\frac
{1}{\sigma^{d}} \biggl(\log\frac{K_2}{\eps\sigma^{d}} \biggr)^{1+d}
\]
for $\eps\in(0,\eps_0)$ and $\sigma\in(0,\sigma_0)$,
with constants $K_1, K_2 >0$ that do not depend on $\eps$ or $\sigma$.
For $\sigma> \sigma_0$, it holds that
\[
\log N(\eps,\GG_{\sigma},\|\cdot\|_{\infty}) \le K_3 \biggl(\log
\frac{1}{\eps} \biggr)^{1+d}
\]
for all $\eps\in(0,\eps_0)$, with $K_3 > 0$ a constant independent
of $\eps$ and $\sigma$.
\end{lem}
\begin{pf}
The first statement is well known; see, for instance, Theo-\break rem~2.7.1~of~\cite{vdVW}.
The second statement is similar to the classical result given by
Theorem~23 of~\cite{Kolmogorov}, which gives the entropy for the class
of analytic functions bounded by a constant on a strip in complex space.
However, the proof of the present statement requires extra care to
identify the role of $\sigma$, because it should not be considered as
an irrelevant constant in our framework.
We omit the details, since the proof of Lemma 4.5 of
\cite{VaartZanten09} is very similar.
\end{pf}

In view of the observations preceding Lemma~\ref{lem: metric entropy},
we now have entropy bounds
for the unit ball of the RKHS in all cases.
Using the results from~\cite{KuelbsLi} and~\cite{LiLinde}, these
translate into results
on the centered small ball probability of $W^{m,\sigma}$.
The first statement of the following lemma follows from the preceding
lemma in combination with the results
of~\cite{LiLinde}.
The second statement is derived from Lemma~\ref{lem: metric entropy}
by arguing as in
the proof of Lemma 4.6 in~\cite{VaartZanten09}.
\begin{lem}\label{bound ball part}
If $d/2 < \gamma< \infty$,
\[
-\log\PP(\|W^{m,\sigma}\|_{\infty}<\eps)\le
K_0 \biggl(\frac{1}{\eps\sigma^{d+\gamma}} \biggr)^{{2d}/({2\gamma-d})}
\]
for all $\eps, \sigma> 0$, with $K_0$ a constant independent of $\eps
$ and $\sigma$.\vadjust{\goodbreak}

If $\gamma= \infty$, there exist $\eps_0, \sigma_0, K_4 > 0$, not
depending on $\eps$ and $\sigma$, such that
\[
-\log\PP(\|W^{m,\sigma}\|_{\infty}<\eps)\le K_4\frac{1}{\sigma
^{d}} \biggl(\log\frac{1}{\eps\sigma^{1+d}} \biggr)^{1+d}
\]
for all $\eps\in(0, \eps_0)$ and $\sigma\in(0,\sigma_0)$.
For $\sigma\ge\sigma_0$ we have
\[
-\log\PP(\|W^{m,\sigma}\|_{\infty}<\eps)\le K_5 \biggl(\log\frac
{1}{\eps} \biggr)^{1+d}
\]
for all $\eps\in(0,\eps_0)$, where $K_5 > 0$ is independent of $\eps
$ and $\sigma$.
\end{lem}

With condition (\ref{eq: toshow3}) in mind, we now consider the
noncentered small ball probabilities of the
process $W^{m,\sigma}$. According to Lemma 5.3 of~\cite{vdVvZRKHS},
we have
for $w_{0} \in C([0,1]^d)$
the inequality
%
%
\begin{equation}\label{eq: noncent}
-\log\PP(\|W^{m,\sigma}-w_{0}\|_{\infty}< 2\eps) \le{\phi
^{m,\sigma}_{w_{0}}(\eps)},
\end{equation}
with $\phi^{m,\sigma}_{w_{0}}$ the so-called concentration function,
defined as follows:
%
%
\begin{equation}\label{concentration function}
\phi^{m,\sigma}_{w_{0}}(\eps)=\inf_{h\in\HHH^{m,\sigma} \dvtx\|
h-\theta_{0}\|_{\infty}\le\eps} \|h\|_{\HHH^{m,\sigma}}^{2}
- \log\PP(\|W^{m,\sigma}\|_{\infty}<\eps).
\end{equation}
(Our function $w_0$ is actually defined only on $\XX$, but we will
extend it to all of $[0,1]^d$ in an
appropriate way later.)
That is to say, the exponent of the noncentered small ball probability
involves the exponent of the centered small ball probability that we
considered above and an approximation term
that quantifies how well $w_0$ can be approximated by elements of the RKHS.

To obtain a suitable approximation, we need an auxiliary result
concerning the approximation of a smooth function $f$ by convolutions.
Define $m_k = \int y^k p(y) \,dy$ for $k \in\NN_0^d$. Next, for $n \in
\NN_0^d$ we
recursively define two
collections of numbers $c_n$ and $d_n$ as follows. If $n. = 1$, we put
$c_n = 0$ and
$d_n = -m_n/n!$. For $n. \ge2$, we define
%
%
\begin{equation}\label{eq: ss}
c_n = -\mathop{\sum_{n = l+k}}_{l. \ge1, k. \ge1} \frac
{(-1)^{k.}}{k!}m_kd_l,\qquad
d_n = \frac{(-1)^{n.} m_n}{n!} + c_n.
\end{equation}
Note that the numbers $c_n$ and $d_n$ are
well defined and that they only depend on the moments of $p$. For a
function $f \in C^\alpha(\RR^d)$ and $\sigma> 0$,
we define the transform $T_{\alpha, \sigma} f$ as follows:
%
%
\begin{equation}\label{eq: sss}
T_{\alpha, \sigma} f = f - \sum_{j=1}^{\beta}\sum_{k.=j}
d_{k}\sigma^{j}(D_{k}^{j}f).
\end{equation}
Here, $\beta$ is the largest integer strictly smaller than $\alpha$
and for a positive integer $j$ and a multi-index $k \in\NN_0^d$ with
$k. = j$,
$D^j_k$ is the $j$th order differential operator
\[
D^j_k = \frac{\partial^j}{\partial x^{k_1}_1 \cdots\partial x^{k_d}_d}.
\]

Let $p_\sigma(x) = \sigma^{-d}p(x/\sigma)$.
\begin{lem}\label{lem: judith}
For $\alpha, \sigma> 0$ and $f \in C^\alpha(\RR^d)$, we have
\[
\|p_\sigma* (T_{\alpha, \sigma} f) - f\|_\infty\le K_6\sigma
^\alpha,
\]
where $K_6> 0 $ is a constant independent of $\sigma$.
\end{lem}

The lemma is an extension of an idea of~\cite{Rou08}, where a similar
method is employed
to approximate arbitrary smooth densities by beta mixtures.
The proof follows the same lines but is somewhat more involved in the
present higher-dimensional case;
see \hyperref[app]{Appendix}.

The following lemma deals with the approximation of the function $w_0$
by elements of the RKHS
of the process $W^{m,\sigma}$.
\begin{lem}
For all $\sigma> 0$, $m \ge1$ and $w_0 \in C^\alpha(\XX)$ there
exists an $h \in\HHH^{m,\sigma}$ such that
$\|h\|_{\HHH^{m,\sigma}} \le K_7 (1\vee\sigma)$ and
\[
\sup_{x \in\XX}|h(x)-w_0(x)| \le\frac{K_8(1\vee\sigma^{\beta
+1})}{\sigma^{1+d}m^{\alpha-\beta}} + K_9\sigma^\alpha,
\]
for $K_7, K_8, K_9 > 0$ constants independent of $\sigma$ and $m$ and
$\beta$ the largest integer
strictly smaller than $\alpha$.
\end{lem}
\begin{pf}
Since $\XX= [a,b]^d \subset(0,1)^d$, we can extend $w_0$ to all of
$\RR^d$ in such a way that
that the resulting function belongs to $C^\alpha(\RR^d)$ and has
support strictly inside $(0,1)^d$.
Using the operator $T_{\alpha, \sigma}$ introduced above [see (\ref
{eq: sss})], we define
\[
h(x) = \sum_{k \in\{1, \ldots, m\}^d} (T_{\alpha, \sigma}w_0)(k/m)
\frac1{m^d}\frac1{\sigma^d}
p \biggl(\frac{x-k/m}{\sigma} \biggr)
\]
for $x \in[0,1]^d$.
By Lemma~\ref{lem: RKHS}, it holds that $h \in\HHH^{m,\sigma}$ and
\[
\|h\|^2_{\HHH^{m,\sigma}} \le\frac1{m^d} \sum_{k \in\{1, \ldots,
m\}^d}
\bigl((T_{\alpha, \sigma} w_0)(k/m) \bigr)^2 \le
\|T_{\alpha, \sigma}w_0\|^2_\infty.
\]
It follows from the definition of $T_{\alpha, \sigma}$ that this
bounded by a constant
times $(1\vee\sigma^\beta)^2$.

It remains to prove the bound for
the approximation error.
By the triangle inequality,
%
%
\begin{equation}\label{eq: karel}
\|h - w_0\|_\infty\le\|h-p_\sigma*(T_{\alpha, \sigma}w_0)\|_\infty
+ \|p_\sigma*(T_{\alpha, \sigma}w_0) - w_0\|_\infty.
\end{equation}
The first term on the right is the difference
between the convolution $p_\sigma*T_{\alpha, \sigma}w_0$ and the
corresponding Riemann sum.
Using again the triangle inequality, we get
\begin{eqnarray*}
&&
|h(x) -(p_\sigma*T_{\alpha, \sigma}w_0)(x)| \\
&&\qquad \le{\sup_{\|y-z\|_\infty\le1/m}} |T_{\alpha, \sigma
}w_0(y)p_\sigma(x-y)-T_{\alpha, \sigma}w_0(z)p_\sigma(x-z)| \\
&&\qquad \le{\|T_{\alpha,
\sigma}w_0\|_\infty\sup_{\|y-z\|_\infty\le1/m}}
|p_\sigma(x-y)-p_\sigma(x-z)|\\
&&\qquad\quad{} +
{\|p_\sigma\|_\infty\sup_{\|y-z\|_\infty\le1/m} }|T_{\alpha, \sigma
}w_0(y)-T_{\alpha, \sigma}w_0(z)|.
\end{eqnarray*}
Now use the facts that $T_{\alpha, \sigma}w_0$ is bounded by a
constant times $1\vee\sigma^\beta$, $p_\sigma$ is bounded
by $\sigma^{-d}$ times a constant, $p$ is Lipschitz and the definition
of $T_{\alpha, \sigma}w_0$
to see that
\[
\|h-p_\sigma*T_{\alpha, \sigma}w_0\|_\infty\le\frac{C_1(1\vee
\sigma^\beta)}{\sigma^{1+d}m}
+ \frac{C_2(1\vee\sigma^\beta)}{\sigma^d m^{\alpha-\beta}}
\le\frac{C_3(1\vee\sigma^{\beta+1})}{\sigma^{1+d}m^{\alpha-\beta}},
\]
which covers the first term on the right-hand side of (\ref{eq: karel}).
Lemma~\ref{lem: judith} implies that the second term is bounded by a
constant times $\sigma^\alpha$.
\end{pf}

By combining the preceding lemma with Lemma~\ref{bound ball part} and
(\ref{eq: noncent}), we obtain the following result.
\begin{lem}\label{lem: noncent}
Let $w_0 \in C^\alpha(\XX)$.

If $\gamma< \infty$,
there exist constants $\eps_0, \sigma_0, K_1, K_2, K_3, K_4 > 0 $,
independent of $\sigma$ and $m$, such that
\[
-\log\PP\Bigl({\sup_{x \in\XX}}|W^{m,\sigma}(x) - w_0(x)| < 2\eps
\Bigr) \le K_1
+K_2 \biggl(\frac{1}{\eps\sigma^{d+\gamma}} \biggr)^{
{2d}/({2\gamma-d})},
\]
provided that
\[
\frac{K_3}{\sigma^{1+d}m^{\alpha- \beta}} + K_4\sigma^\alpha<
\eps< \eps_0
\]
and $\sigma\in(0,\sigma_0)$.

If $\gamma= \infty$, there exist constants $\eps_0, \sigma_0, K_1,
K_2, K_3, K_4 > 0 $, independent of $\sigma$ and $m$, such that
\[
-\log\PP\Bigl({\sup_{x \in\XX}}|W^{m,\sigma}(x) - w_0(x)| < 2\eps
\Bigr) \le K_1
+ K_2\frac{1}{\sigma^{d}} \biggl(\log\frac{1}{\eps\sigma^{1+d}}
\biggr)^{1+d},
\]
provided that
\[
\frac{K_3}{\sigma^{1+d}m^{\alpha- \beta}} + K_4\sigma^\alpha<
\eps< \eps_0
\]
and $\sigma\in(0,\sigma_0)$.
\end{lem}

\subsection{\texorpdfstring{Proof of Theorem \protect\ref{thm: general}}{Proof of Theorem 2.2}}

\subsubsection{\texorpdfstring{Condition (\protect\ref{eq: toshow3})}{Condition (2.5)}}
By definition of the process $W$ and conditioning,
\begin{eqnarray*}
&& \PP\Bigl({\sup_{x \in\XX}} |W(x)-w_0(x)| \le\eps\Bigr)\\
&&\qquad = \sum_{m =1}^\infty\lambda_m \int_0^\infty g(\sigma)
\PP\Bigl({\sup_{x \in\XX} }|W^{m,\sigma}(x) - w_0(x)| < \eps\Bigr)\,
d\sigma,
\end{eqnarray*}
where $\lambda_m = \PP(M= m)$. If $\gamma< \infty$, Lemma \ref
{lem: noncent} implies that there exist constants $\eps_0, C_1, C_2,
C_3, C_4 > 0$, independent
of $\sigma$ and $m$, such that
if $\eps< \eps_0$ and
\[
\tfrac12 C_1\eps^{1/\alpha} < \sigma< C_1\eps^{1/\alpha} \le1,\qquad
m \ge C_2\eps^{-({1+d+\alpha})/({\alpha(\alpha-\beta)})},
\]
then
\[
-\log\PP\Bigl({\sup_{x \in\XX} }|W^{m,\sigma}(x) - w_0(x)| < \eps
\Bigr) \le C_3
+ C_4 \biggl(\frac{1}{\eps\sigma^{d+\gamma}} \biggr)^{{2d}/({2\gamma-d})}.
\]
Hence, the probability of interest is bounded from below, for $\eps<
\eps_0$, by
\begin{eqnarray*}
&& e^{-C_3} \sum_{m \ge C_2\eps^{-({1+d+\alpha
})/({\alpha(\alpha-\beta)})}}\lambda_m
\int_{C_1\eps^{1/\alpha}/2}^{C_1\eps^{1/\alpha}}g(\sigma)\exp
\biggl(-C_4 \biggl(\frac{1}{\eps\sigma^{d+\gamma}} \biggr)^{
{2d}/({2\gamma-d})} \biggr) \,d\sigma\\
&&\qquad \ge C_5 \exp
\bigl(-C_6\eps^{-({\alpha+d+\gamma})/{\alpha}{2d}/({2\gamma
-d})} \bigr)
\end{eqnarray*}
for constants $C_5, C_6 > 0$.
It follows that condition (\ref{eq: toshow3}) is fulfilled for
%
%
\begin{equation}\label{eq: epsnbis}
\eps_n = M_1{n}^{-{\alpha}/({d_\gamma+2\alpha(1+\delta_\gamma)})}
\end{equation}
for $M_1 > 0$ an appropriate constant and $d_\gamma= 2d(d+\gamma
)/(2\gamma-d)$,
$\delta_\gamma= d/(2\gamma- d)$.

If $\gamma= \infty$,
the same reasoning implies that
there exist constants $C_5, C_6 > 0$ such that, for $\eps> 0$ small enough,
\[
\PP\Bigl({\sup_{x \in\XX} }|W(x)-w_0(x)| \le\eps\Bigr) \ge C_5
e^{-C_6 \eps^{-d/\alpha}\log^{r\vee(1+d)}(1/\eps)}.
\]
It follows that, in this case, condition (\ref{eq: toshow3}) is
fulfilled for
%
%
\begin{equation}\label{eq: epsn}
\eps_n = M_1{n}^{-{\alpha}/({d+2\alpha})} \log^t n
\end{equation}
for $M_1 > 0$ an appropriate constant, provided that $t \ge(r\vee
(1+d))/(2+d/\alpha)$.

\subsubsection{\texorpdfstring{Construction of the sets $B_n$ and condition (\protect\ref{eq: toshow2})}
{Construction of the sets $B_n$ and condition (2.4)}}

First, suppose that $\gamma< \infty$ again.
For $L, R, \eps> 0$, we define
\[
B = LC^\gamma_{{R^{-(d+\gamma)}}{\|p\|_\gamma}}([0,1]^d) + \eps\BBB_1,\vadjust{\goodbreak}
\]
where $\BBB_1$ is the unit ball of the space $C([0,1]^d)$.
The sieves $B_n$ will be defined by making appropriate choices for the
$L,R$ and $\eps$ below.
Recall that in this case
$\HHH^{m,\sigma}_1 \subset C^\gamma_{{\sigma^{-(d+\gamma)}}{\|p\|
_\gamma}}([0,1]^d)$.
Hence, by the Borell--Sudakov inequality (see, e.g.,
\cite{Lifshits}), with $\Phi$ the standard normal distribution function and
for $\sigma\ge R$,
\begin{eqnarray*}
\PP(W^{m,\sigma} \notin B) & \le & \PP(W^{m,\sigma} \notin L\HHH
^{m,\sigma}_1 + \eps\BBB_1)\\
& \le & 1-\Phi\bigl(\Phi^{-1}\bigl(\PP(\|W^{m,\sigma}\|_\infty\le\eps)\bigr) + L\bigr).
\end{eqnarray*}
By Lemma~\ref{bound ball part}, we have, for
$\sigma\ge R$ and $R \le1$,
\[
\PP(\|W^{m,\sigma}\|_\infty\le\eps) \ge e^{-K_6 R^{-d_\gamma}\eps
^{-2d/(2\gamma- d)}}
\]
for a constant $K_6 > 0$ and $\eps>0$ small enough.
Since $\Phi^{-1}(y) \ge-\sqrt{(5/2) \log(1/y)}$ for $y \in(0,1/2)$,
it follows that
\begin{eqnarray*}
\PP(W^{m,\sigma} \notin B) & \le & 1-\Phi\bigl(L-\sqrt{(5/2)K_6
R^{-d_\gamma}\eps^{-2d/(2\gamma- d)}} \bigr) \\
& \le & e^{-1/2(L-\sqrt{(5/2)K_6 R^{-d_\gamma}\eps^{-2d/(2\gamma
- d)}})^2},
\end{eqnarray*}
for $\sigma\ge R$ and $L \ge\sqrt{(5/2)K_6 R^{-d_\gamma}\eps
^{-2d/(2\gamma- d)}}$.
By the definition of $W$ and conditioning,
\[
\PP(W \notin B) \le\sum_{m=1}^\infty\lambda_m \int_{R}^\infty
g(\sigma)\PP(W^{m,\sigma} \notin B) \,d\sigma
+ \PP(\Sigma< R).
\]
By the preceding, the first term on the right is bounded by
\[
e^{-1/2(L-\sqrt{(5/2)K_6 R^{-d_\gamma}\eps^{-2d/(2\gamma- d)}})^2}.
\]
The assumption on $g$ and a
substitution show that the second term is bounded by
\[
D_3\int_{1/R}^\infty x^{q-2}e^{-D_4x^{d_\gamma}(\log x)^{r}} \,dx.
\]
By Lemma 4.9 of~\cite{VaartZanten09}, this is further bounded by
\[
\frac{2D_3}{dD_4}\frac{(1/R)^{q-2-d_\gamma+1}}{(\log
(1/R))^{r}}e^{-D_4(1/R)^{d_\gamma}(\log(1/R))^{r}}
\le e^{-1/2D_4(1/R)^{d_\gamma}(\log(1/R))^r}
\]
for $R$ small enough.

Given $C > 1$, we now define the sieve $B_n$ by
\[
B_n = L_nC^\gamma_{{R_n^{-(d+\gamma)}}{\|p\|_\gamma}}([0,1]^d) +
\eps_n\BBB_1,
\]
where $\eps_n$ is given by (\ref{eq: epsnbis}). To show that (\ref
{eq: toshow2}) holds, we have to show we can choose $R_n$ and $L_n$
such that
\[
\frac1{R_n^{d_\gamma}}\log^r \frac1{R_n} \ge Cn\eps^2_n
\]
and
\[
\bigl(L_n-\sqrt{(5/2)K_6 R_n^{-d_\gamma}\eps_n^{-2d/(2\gamma- d)}}\bigr)^2
\ge Cn\eps^2_n.
\]
Observe that if we take
\[
\frac1{R_n^{d_\gamma}} = Mn^{({d_\gamma+2\alpha\delta_\gamma
})/({d_\gamma+2\alpha(1+\delta_\gamma)})}
\]
for a large enough constant $M$, the first condition is satisfied.
The second condition is then fulfilled if we choose
\[
L^2_n = N n^{({d_\gamma+ 4\alpha\delta_\gamma})/({d_\gamma
+2\alpha(1+\delta_\gamma)})}
\]
for $N$ large enough.

Next, we consider the case $\gamma= \infty$.
Recall that $\GG_\sigma$ is the set of all analytic functions defined
on the strip
$S_{\sigma}=\{z\in\CCC^{d}\dvtx\forall j $ $ |\Im z_{j}|\le\sigma\}$
that are bounded by
$K\sigma^{-d}$ on $S_\sigma$. Arguing as before and now using that
$\HHH^{m,\sigma}_1 \subset\GG_\sigma$ and
$\GG_{\sigma_1} \sub\GG_{\sigma_2}$ if
$\sigma_1 \ge\sigma_2$, we get,
for $L, R, \eps> 0$
and $B = L\GG_{R} + \eps\BBB_1$,
\[
\PP(W^{m,\sigma} \notin B) \le e^{-1/2(L-\sqrt{(5/2)K_6
R^{-d}(\log(1/(\eps R^{1+d})))^{1+d}})^2}
\]
for $\sigma\ge R$ and $L \ge\sqrt{(5/2)K_6 R^{-d}(\log(1/(\eps
R^{1+d})))^{1+d}}$.
By the same conditioning argument as before, it follows that if, given
$C > 1$, we define $B_n$ in this case by
\[
B_n = L_n\GG_{R_n} +\eps_n\BBB_1,
\]
where $\eps_n$ is given by (\ref{eq: epsn}),
then condition (\ref{eq: toshow2}) is fulfilled
if we choose $R_n$ and $L_n$ such that
\[
\frac1{R_n^d}\log^r \frac1{R_n} \ge Cn\eps^2_n
\]
and
\[
\bigl(L_n-\sqrt{(5/2)K_6 R_n^{-d}\bigl(\log\bigl(1/(\eps_n
R_n^{1+d})\bigr)\bigr)^{1+d}} \bigr)^2 \ge Cn\eps^2_n.
\]
Observe that we can take
\[
\frac1{R_n^d} = Mn^{{d}/({d+2\alpha})}\log^v n
\]
for a large enough constant $M$ and $v \ge2t-r$ [with $t$ as in (\ref
{eq: epsn})], and $L_n$ a large enough power of $n$.

\subsubsection{Entropy condition}

Suppose $\gamma< \infty$.
For the entropy of the sieve $B_n$, we have in this case, for
$\overline{\eps}_n \ge\eps_n$,
\begin{eqnarray*}
N(2\overline\eps_n, B_n, \|\cdot\|_\infty) & \le & N\bigl(\overline\eps
_n, L_nC^\gamma_{{R_n^{-(d+\gamma)}}{\|p\|_\gamma}}([0,1]^d), \|
\cdot\|_\infty\bigr)\\
& \le & N\bigl(\overline\eps_n R_n^{d+\gamma}/({L_n}{\|p\|_\gamma}),
C^\gamma_{{1}}([0,1]^d), \|\cdot\|_\infty\bigr).
\end{eqnarray*}
Hence  (see Lemma~\ref{lem: metric entropy}),
\[
\log N(2\overline\eps_n, B_n, \|\cdot\|_\infty) \le
K_1 \biggl( \frac{L_n}{\overline\eps_n R_n^{d+\gamma}}
\biggr)^{d/\gamma}.
\]
This is bounded by a constant times $n\overline\eps_n^2$ for
\[
\overline\eps_n \gtrsim\frac{L_n^{d/(d+2\gamma)}}{n^{
{\gamma}/({d+2\gamma})}R_n^{{d(d+\gamma)}/({d+2\gamma})}}.
\]
For $L_n$ and $R_n$ chosen as above, this yields
\[
\overline\eps_n \gtrsim n^{-\frac{\alpha(1-(d\delta_\gamma
)/(2\gamma))}
{d_\gamma+ 2\alpha(1+\delta_\gamma) + d(d_\gamma+ 2\alpha
(1+\delta_\gamma))/(2\gamma)}}.
\]
Note that $\overline\eps_n$ is always larger
than $\eps_n$, as was required.

Let now $\gamma= \infty$. Arguing as before, we have in this case,
for $\overline{\eps}_n \ge\eps_n$,
\[
N(2\overline\eps_n, B_n, \|\cdot\|_\infty) \le
N(\overline\eps_n/{L_n}, \GG_{R_n}, \|\cdot\|_\infty)
\le
K_1\frac{1}{R_n^{d}} \biggl(\log\frac{L_n}{\overline\eps_nR_n^{d}}
\biggr)^{1+d}
\]
by Lemma~\ref{lem: metric entropy}.
With the choices of $R_n$ and $L_n$ made in this case above
and for $\overline\eps_n$ bounded from below by a power of $n$, this
is bounded by a constant times
$n^{{d}/({d+2\alpha})}\log^{1+d+v} n$.
This is further bounded by a constant times $n\overline\eps_n^2$ for
\[
\overline\eps_n = n^{-{\alpha}/({d+2\alpha})}\log^{a} n,
\]
provided $a \ge(1+d+v)/2$. The requirement that $\overline{\eps}_n
\ge\eps_n$
translates into the condition $a \ge t$.

\begin{appendix}\label{app}
\section*{Appendix}

\begin{pf*}{Proof of Lemma~\ref{lem: judith}}
The proof is by induction on $\beta$, which is the largest integer
strictly smaller than $\alpha$.
If $\beta= 0$ then $\alpha\in(0,1]$ and $T_{\alpha, \sigma} f = f$
and the
statement of the claim is standard.
To prove the induction step, suppose now that $\beta\ge1$. By
definition of
$T_{\alpha, \sigma} f$, we have
\begin{eqnarray*}
&& (p_\sigma*T_{\alpha, \sigma} f- f)(x)\\
&&\qquad = \int p_\sigma(y) \Biggl(f(x-y)-f(x) -
\sum_{j=1}^{\beta}\sum_{k.=j} d_{k}\sigma^{j}(D_{k}^{j}f)(x-y)\Biggr)\, dy.
\end{eqnarray*}
By Taylor's formula and the fact that $f\in C^\alpha$,
\[
f(x-y)-f(x) = \sum_{j=1}^{\beta}\sum_{k.=j} \frac
{(-y)^k}{k!}(D_{k}^{j}f)(x) + R(x,y),
\]
where $|R(x,y)| \le C \|y\|^\alpha$. It follows that
\begin{eqnarray*}
&&(p_\sigma*T_{\alpha, \sigma} f - f)(x) \\[-2pt]
&&\qquad=  \int p_\sigma(y) R(x,y) \,dy\\[-2pt]
&&\qquad\quad{} + \sum_{j=1}^{\beta}\sum_{k.=j}
\biggl( \frac{1}{k!}(-1)^j (D^j_kf)(x)\sigma^j m_k - d_k\sigma^j
\bigl(p_\sigma* (D_k^jf) \bigr)(x) \biggr).
\end{eqnarray*}
The first term on the right is easily seen to be bounded by a constant
times $\sigma^\alpha$. To
see that this holds for the second term as well, we use the induction
hypothesis.

By definition of the constants $c_k$ and $d_k$ [see (\ref{eq: ss})],
the second term can be written as
\[
\sum_{j=1}^{\beta}\sum_{k_{\cdot}=j}
\biggl( \frac{(-1)^j}{k!} \sigma^j m_k \bigl(D^j_kf - p_\sigma*
(D_k^jf) \bigr)(x) - c_k\sigma^j
\bigl(p_\sigma* (D_k^jf) \bigr)(x) \biggr).
\]
Now for $j \le\beta$ and $k. = j$, consider the decomposition
\begin{eqnarray*}
&& D^j_kf - p_\sigma* (D_k^jf) \\[-2pt]
&&\qquad = \bigl(D^j_kf - p_\sigma* (T_{\alpha-j, \sigma}{D_k^jf})
\bigr)\\[-2pt]
&&\quad \qquad{}+ \bigl(p_\sigma* (T_{\alpha-j, \sigma}{D_k^jf}) - p_\sigma*
(D_k^jf) \bigr).
\end{eqnarray*}
Since $D^j_kf \in C^{\alpha- j}$, the induction hypothesis implies
that the first
term on the right is uniformly bounded by a constant times $\sigma
^{\alpha-j}$.
Combined with the first\ display of the paragraph, this shows that it
suffices to show that
\[
\sum_{j=1}^{\beta}\sum_{k.=j}
\biggl( \frac{(-1)^j}{k!} \sigma^j m_k (T_{\alpha-j, \sigma
}{D_k^jf} - D_k^jf ) - c_k\sigma^j
(D_k^jf ) \biggr) = 0
\]
identically.
Straightforward algebra shows that
\[
T_{\alpha-j, \sigma}{D_k^jf} - D_k^jf = - \sum_{i=1}^{\beta-j}\sum
_{l.=i} d_{l}\sigma^{i}D^{i+j}_{k+l}f.
\]
Hence,
\begin{eqnarray*}
&&\sum_{j=1}^{\beta}\sum_{k.=j}
\frac{(-1)^j}{k!} \sigma^j m_k (T_{\alpha-j, \sigma}{D_k^jf}
- D_k^jf ) \\[-2pt]
&&\qquad = - \sum_{j=1}^{\beta}\sum_{k.=j}\sum_{i=1}^{{\beta- j}}\sum_{l.=i}
\frac{(-1)^j}{k!} m_k d_{l}\sigma^{i+j}D_{l+k}^{i+j}f\\[-2pt]
&&\qquad = - \sum_{s = 2}^\beta\sum_{n. = s}
\biggl(\mathop{\sum_{n=l+k}}_{l. \ge1, k. \ge1} \frac{(-1)^{k.}}{k!} m_k
d_{l} \biggr)\sigma^s D^s_nf.
\end{eqnarray*}
By definition of the numbers $c_n$ and $d_n$ this equals
\[
\sum_{s=1}^{\beta}\sum_{n.=s}
c_n\sigma^s D_n^sf,
\]
and the proof is complete.
\end{pf*}
\end{appendix}

\printaddresses

\end{document}